\documentclass{article}
\usepackage{amsfonts}
\usepackage{amsmath}

\newtheorem{theorem}{Theorem}

\newtheorem{corollary}[theorem]{Corollary}

\newtheorem{lemma}[theorem]{Lemma}

\begin{document}

\title{Path to survival for the critical branching processes in a random
environment\thanks{%
This work is supported by the RFBR under the grant N 14-01-00318. } }
\author{Vatutin V.A.\thanks{%
Department of Discrete Mathematics, Steklov Mathematical Institute, 8,
Gubkin str., 119991, Moscow, Russia; e-mail: vatutin@mi.ras.ru}, Dyakonova
E.E.\thanks{%
Department of Discrete Mathematics, Steklov Mathematical Institute, 8,
Gubkin str., 119991, Moscow, Russia; e-mail: elena@mi.ras.ru}}
\date{}
\maketitle

\begin{abstract}
A critical branching process $\left\{ Z_{k},k=0,1,2,...\right\} $ in a
random environment is considered.\ A conditional functional limit theorem
for the properly scaled process $\left\{ \log Z_{pu},0\leq u<\infty \right\}
$ is established\ under the assumptions $Z_{n}>0$ and $p\ll n$. It is shown
that the limiting process is a Levy process conditioned to stay nonnegative.
The proof of this result is based on a limit theorem describing the
distribution of the initial part of the trajectories of a driftless random
walk conditioned to stay nonnegative.

\medskip


\end{abstract}

\section{Introduction}

We consider a branching process in a random environment specified by a
sequence of independent identically distributed random laws. Denote by $%
\Delta $ the space of probability measures on $\mathbb{N}_{0}=\left\{
0,1,2,...\right\} $. Equipped with the metric of total variation, $\Delta $
becomes a Polish space. Let $Q$ be a random variable taking values in $%
\Delta $. Then, an infinite sequence
\begin{equation}
\Pi =(Q_{1},Q_{2},\ldots )  \label{DefEnvir}
\end{equation}%
of i.i.d. copies of $Q$ is said to form a \emph{random environment}. A
sequence of $\mathbb{N}_{0}$-valued random variables $Z_{0},Z_{1},\ldots $
is called a \emph{branching process in the random environment} $\Pi $, if $%
Z_{0}$ is independent of $\Pi $ and given $\Pi $ the process $\mathcal{Z}%
=(Z_{0},Z_{1},\ldots )$ is a Markov chain with
\begin{equation}
\mathcal{L}\left( Z_{n}\;|\;Z_{n-1}=z,\,\Pi =(q_{1},q_{2},\ldots )\right) \
=\ \mathcal{L}\left( \xi _{n1}+\cdots +\xi _{nz}\right)   \label{transition}
\end{equation}%
for every $n\geq 1,\,z\in \mathbb{N}_{0}$ and $q_{1},q_{2},\ldots \in \Delta
$, where $\xi _{n1},\xi _{n2},\ldots $ are i.i.d. random variables with
distribution $q_{n}$.

In the language of branching processes $Z_{n}$ is the $n$th generation size
of the population and $Q_{n}$ is the distribution of the number of children
of an individual at generation $n-1$. We assume that $Z_{0}=1$ a.s. for
convenience and denote the corresponding probability measure on the
underlying probability space by~$\mathbf{P}$. (If we refer to other
probability spaces, then we use notation $\mathbb{P}$, $\mathbb{E}$ \ and $%
\mathbb{L}$ for the respective probability measures, expectations and laws.)

As it turns out the properties of $\mathcal{Z}$ are first of all determined
by its associated random walk $\mathcal{S}:=\left\{ S_{n},n\geq 0\right\} $.
This random walk has initial state $S_{0}=0$ and increments $%
X_{n}=S_{n}-S_{n-1}$, $n\geq 1$ defined as
\begin{equation*}
X_{n}\ :=\ \log \ \left( \sum_{y=0}^{\infty }y\ Q_{n}(\{y\})\right) ,
\end{equation*}%
which are i.i.d. copies of the logarithmic mean offspring number
\begin{equation*}
X\ :=\ \log \ \left( \sum_{y=0}^{\infty }y\ Q(\{y\})\right) .
\end{equation*}

Following \cite{AGKV05} we call the process $\mathcal{Z}:=\left\{
Z_{n},\,n\geq 0\right\} $ \textit{critical} if and only if the random walk $%
\mathcal{S}$ is oscillating, that is,
\begin{equation*}
\limsup_{n\rightarrow \infty }S_{n}=\infty \ \text{ and }\
\liminf_{n\rightarrow \infty }S_{n}=-\infty .
\end{equation*}%
It is shown in \cite{AGKV05} that the extinction moment of the critical
branching process in a random environment is finite with probability $1$.
For this reason it is natural to study the asymptotic behavior of the
survival probability $\mathbf{P}(Z_{n}>0)$ as $n\rightarrow \infty .$ This
has been done in \cite{AGKV05}: If
\begin{equation}
\lim_{n\rightarrow \infty }\mathbf{P}\left( S_{n}>0\right) =\rho \in (0,1),
\label{Spit}
\end{equation}%
then (under some mild additional assumptions to be specified later on)
\begin{equation}
\mathbf{P}(Z_{n}>0)\sim \theta \mathbf{P}(\min \left(
S_{0},S_{1},...,S_{n}\right) \geq 0)=\theta \frac{l(n)}{n^{1-\rho }},
\label{tails}
\end{equation}%
where $l(n)$ is a slowly varying function and $\theta $ is a known positive
constant whose explicit expression is given by formula (\ref{DefTheta})
below.

Let
\begin{equation*}
\mathcal{A}=\{0<\alpha <1;\,\left\vert \beta \right\vert <1\}\cup \{1<\alpha
<2;\left\vert \beta \right\vert <1\}\cup \{\alpha =1,\beta =0\}\cup \{\alpha
=2,\beta =0\}
\end{equation*}%
be a subset in $\mathbb{R}^{2}.$ For $(\alpha ,\beta )\in \mathcal{A}$ and a
random variable $X$ write $X\in \mathcal{D}\left( \alpha ,\beta \right) $ if
the distribution of $X$ belongs to the domain of attraction of a stable law
with characteristic function%
\begin{equation}
\mathcal{H}_{\alpha ,\beta }\mathbb{(}t\mathbb{)}:=\exp \left\{
-c|t|^{\,\alpha }\left( 1+i\beta \frac{t}{|t|}\tan \frac{\pi \alpha }{2}%
\right) \right\} ,\ c>0,  \label{std}
\end{equation}%
and, in addition, $\mathbf{E}X=0$ if this moment exists.

\ Denote $\mathbb{N}_{+}:=\left\{ 1,2,...\right\} $ and let $\left\{
c_{n},n\geq 1\right\} $ be a sequence of positive integers specified by the
relation%
\begin{equation}
c_{n}:=\inf \left\{ u\geq 0:G(u)\leq n^{-1}\right\} ,  \label{Defa}
\end{equation}%
where
\begin{equation*}
G(u):=\frac{1}{u^{2}}\int_{-u}^{u}x^{2}\mathbf{P}(X\in dx).
\end{equation*}%
It is known (see, for instance, \cite[Ch. XVII, \S 5]{FE}) that, for every $%
X\in \mathcal{D}(\alpha ,\beta )$ the function $G(u)$ is regularly varying
with index $-\alpha $. This implies that $\left\{ c_{n},n\geq 1\right\} $ is
a regularly varying sequence with index $\alpha ^{-1}$, i.e., there exists a
function $l_{1}(n),$ slowly varying at infinity, such that
\begin{equation}
c_{n}=n^{1/\alpha }l_{1}(n).  \label{asyma}
\end{equation}%
In addition, the scaled sequence $\left\{ S_{n}/c_{n},\,n\geq 1\right\} $
converges in distribution, as $n~\rightarrow ~\infty ,$ to the stable law
given by (\ref{std}). \

Observe that if $X\in \mathcal{D}\left( \alpha ,\beta \right) ,$ then (see,$%
\ $for instance, \cite{Zol57}) the quantity $\rho $ in (\ref{Spit}) is
calculated by the formula
\begin{equation}
\displaystyle\rho =\left\{
\begin{array}{ll}
\frac{1}{2},\ \text{if \ }\alpha =1\text{ or }2 &  \\
\frac{1}{2}+\frac{1}{\pi \alpha }\arctan \left( \beta \tan \frac{\pi \alpha
}{2}\right) ,\text{ otherwise}. &
\end{array}%
\right.  \label{ro}
\end{equation}%
In particular, $\rho \in \left( 0,1\right) $.

Denote%
\begin{equation*}
M_{n}:=\max \left( S_{1},...,S_{n}\right) ,\quad L_{k,n}:=\min_{k\leq j\leq
n}S_{j},\quad L_{n}:=L_{0,n}=\min \left( S_{0},S_{1},...,S_{n}\right)
\end{equation*}%
and introduce a right-continuous renewal function
\begin{equation}
V(x):=1+\sum_{k=1}^{\infty }\mathbf{P}\left( -S_{k}\leq x,M_{k}<0\right) ,\
x\geq 0,  \label{ren}
\end{equation}%
and $0$ elsewhere. In particular, $V(0)=1$.

The fundamental property of $V$ is the identity
\begin{equation}
\mathbf{E}[V(x+X);x+X\geq 0]\ =\ V(x),\text{ }x\geq 0\ ,  \label{harm}
\end{equation}%
which holds for any oscillating random walk.

It follows from (\ref{harm}) that $V$ gives rise to further probability
measures $\mathbf{P}_{x}^{+},x\geq 0,$ specified by corresponding
expectations $\mathbf{E}_{x}^{+}$. The construction procedure of this
measure is explained in \cite{AGKV05} in detail. We only recall that if the
random walk $\mathcal{S}=(S_{n},n\geq 0)$ with $S_{0}=x\geq 0$ is adapted to
some filtration $\mathcal{F}=(\mathcal{F}_{n})$ and $\zeta _{0},\zeta
_{1},\ldots $ is a sequence of random variables, adapted to $\mathcal{F}$,
then for each fixed $n$ and a bounded and measurable function $g_{n}:\mathbb{%
R}^{n+1}\rightarrow \mathbb{R}$,%
\begin{equation}
\mathbf{E}_{x}^{+}[g_{n}(\zeta _{0},\ldots ,\zeta _{n})]\ :=\ \frac{1}{V(x)}%
\mathbf{E}_{x}[g_{n}(\zeta _{0},\ldots ,\zeta _{n})V(S_{n});L_{n}\geq 0],
\notag
\end{equation}%
where $\mathbf{E}_{x}$ is the expectation corresponding to the probability
measure $\mathbf{P}_{x}$ which is generated by $\mathcal{S}$. Under the
measure $\mathbf{P}^{+}=\mathbf{P}_{0}^{+}$ the sequence $S_{0},S_{1},\ldots
$ is a Markov chain with state space $[0,\infty )$ and transition
probabilities
\begin{equation*}
\mathbf{P}^{{+}}(x,dy)\ :=\ \frac{1}{V(x)}\mathbf{P}\left( x+X\in dy\right)
V(y),\quad x\geq 0\ .
\end{equation*}%
It is the random walk conditioned never to enter $(-\infty ,0).$

We now describe in brief a construction of Levy processes conditioned to
stay positive following basically the definitions given in \cite{Chau06} and
\cite{CCh2006}.

Let $\Omega :=D\left( [0,\infty ),\mathbb{R}\right) $ be the space of
real-valued c\`{a}dl\`{a}g paths on the real half-line $[0,\infty )$ and let
$\mathcal{B}:=\left\{ B_{t},t\geq 0\right\} $ be the coordinate process
defined by the equality $B_{t}(\omega )=\omega _{t}$ for $\omega \in \Omega $%
. In the sequel we consider also the spaces $\Omega _{U}:=D\left( [0,U],%
\mathbb{R}\right) ,U>0.$

We endow the spaces $\Omega $ and $\Omega _{U}$ with Skorokhod topology and
denote by $\mathcal{F}=\left\{ \mathcal{F}_{t},t\geq 0\right\} $ and by $%
\mathcal{F}^{U}\mathcal{=}\left\{ \mathcal{F}_{t},t\in \lbrack 0,U]\right\} $
(with some misuse of notation) the natural filtrations of the processes $%
\mathcal{B}$ and $\mathcal{B}^{U}=\left\{ B_{t},t\in \lbrack 0,U]\right\} $.

Let $\mathbb{P}_{x}$ be the law on $\Omega $ an $\alpha -$stable process $%
\mathcal{B}$, $\alpha \in (0,2]$ started at $x$ and let $\mathbb{P}=\mathbb{P%
}_{0}$. Denote by $\rho =\mathbb{P}\left( B_{1}\geq 0\right) $ the
positivity parameter of the process $\mathcal{B}$ (in fact, this quantity is
the same as in (\ref{ro})). We now introduce an analogue of the measure $%
\mathbf{P}^{+}$ for Levy processes. Namely, following \cite{Chau97} we
specify for all $t>0,$ $\mathcal{A}\in \mathcal{F}_{t}$ the law $\mathbb{P}%
_{x}^{+}$ on $\Omega $ \ of the Levy process starting at point $x>0$ and
conditioned to stay positive by the equality
\begin{equation*}
\mathbb{P}_{x}^{+}\left( A\right) :=\frac{1}{x^{\alpha \left( 1-\rho \right)
}}\mathbb{E}_{x}\left[ B_{t}^{\alpha \left( 1-\rho \right) }I\left\{
\mathcal{A}\right\} I\left\{ \inf_{0\leq u\leq t}B_{u}\geq 0\right\} \right]
,
\end{equation*}%
where $I\left\{ \mathcal{C}\right\} $ is the indicator of the event $%
\mathcal{C}$.

Thus, $\mathbb{P}_{x}^{+}$ is an $h-$transform of the Levy process killed
when it first enters the negative half-line. The corresponding positive
invariant function is $H(x)=x^{\alpha \left( 1-\rho \right) }$.

This definition has no sense for $x=0$. However, it is shown in \cite{Chau06}
that it is possible to construct a law $\mathbb{P}^{+}:=\mathbb{P}_{0}^{+}$
and a c\`{a}dl\`{a}g Markov process with the same semigroup as $\left(
\mathcal{B},\left\{ \mathbb{P}_{x}^{+},x>0\right\} \right) $ and such that $%
\mathbb{P}^{+}\left( B_{0}=0\right) =1$. Moreover,%
\begin{equation*}
\mathbb{P}_{x}^{+}\Longrightarrow \mathbb{P}^{+},\text{ as }x\downarrow 0,%
\end{equation*}%
where here and in what follows $\Longrightarrow $ means weak convergence.

\bigskip Let $\mathbb{P}^{(m)}$ be the law on $\Omega _{1}$ of the meander
of length 1 associated with $\left( \mathcal{B},\mathbb{P}\right) ,$ i.e.
\begin{equation}
\mathbb{P}^{(m)}\left( \cdot \right) :=\lim_{x\downarrow 0}\mathbb{P}%
_{x}\left( \cdot \Big|\inf_{0\leq u\leq 1}B_{u}\geq 0\right) .  \label{MesPm}
\end{equation}%
Thus, the law $\mathbb{P}^{(m)}$ may be viewed as the law of the Levy
process $\left( \mathcal{B},\mathbb{P}\right) $ conditioned to stay
nonnegative on the time-interval $\left( 0,1\right) $ while the law $\mathbb{%
P}^{+}$ corresponds to the law of the Levy process conditioned to stay
nonnegative on the whole real half-line $(0,\infty )$.

It is proved in \cite{Chau06} that $\mathbb{P}^{(m)}$ and $\mathbb{P}^{+}$
are absolutely continuous with respect to each other: for every event $%
\mathcal{A}\in \mathcal{F}_{1}$
\begin{equation}
\mathbb{P}^{+}\left( A\right) =C_{0}\mathbb{E}^{(m)}\left[ I\left\{ \mathcal{%
A}\right\} B_{1}^{\alpha \left( 1-\rho \right) }\right] ,  \label{DefPplus}
\end{equation}%
where (see, for instance, formulas (3.5), (3.6), and (3.11) in \cite{CCh2006}%
)
\begin{equation}
C_{0}:=\lim_{n\rightarrow \infty }V(c_{n})\mathbf{P}\left( L_{n}\geq
0\right) \in \left( 0,\infty \right) .  \label{AsH}
\end{equation}%
Hence,%
\begin{equation}
C_{0}^{-1}=\mathbb{E}^{(m)}\left[ B_{1}^{\alpha \left( 1-\rho \right) }%
\right] .  \label{Cexplis}
\end{equation}%
In fact, one may extend the absolute continuity given in (\ref{DefPplus}) to
an arbitrary interval $[0,U]$ be considering the respective space $\Omega
_{U}$ instead of $\Omega _{1}$ and conditioning by the event $\inf_{0\leq
u\leq U}B_{u}\geq 0$ in (\ref{MesPm}).

Set%
\begin{equation*}
\zeta (a):=\frac{\sum_{y=a}^{\infty }y^{2}Q\left( \left\{ y\right\} \right)
}{\left( \sum_{y=0}^{\infty }y\ Q(\{y\})\right) ^{2}},\ a\in \mathbb{N}_{0}.
\end{equation*}

In what follows we say that

1) \textit{Condition }$A1$\textit{\ is valid} if $X\in \mathcal{D}\left(
\alpha ,\beta \right) ;$

2) \textit{Condition }$A2$\textit{\ is valid} if
\begin{equation}
\mathbf{E}\left( \log ^{+}\zeta (a)\right) ^{\alpha +\varepsilon }<\infty
\label{Neeepsil}
\end{equation}%
for some $\varepsilon >0$ and $a\in \mathbb{N}_{0}$;

3) \textit{Condition }$A$\textit{\ is valid} if Conditions $A1$ and $A2$ hold
true and, in addition, the parameter $p=p(n)$ tends to infinity as $%
n\rightarrow \infty $ in such a way that%
\begin{equation}
\lim_{n\rightarrow \infty }n^{-1}p=\lim_{n\rightarrow \infty }n^{-1}p(n)=0.
\label{Msmalln}
\end{equation}

Introduce two processes%
\begin{equation*}
\mathcal{H}^{p}:=\left\{ \frac{\log Z_{\left[ pu\right] }}{c_{p}},0\leq
u<\infty \right\} ,\quad \mathcal{G}^{n}:=\left\{ \frac{\log Z_{\left[ nt%
\right] }}{c_{n}},0\leq t\leq 1\right\} .
\end{equation*}%
We are now ready to formulate two main results of the paper.

The first theorem describes the initial stage of the trajectories of the
critical branching process in a random environment that provide survival of
the process for a long time:

\begin{theorem}
\label{T_small} If Condition $A$ is valid, then, as $n\rightarrow \infty $%
\begin{equation*}
\mathcal{L}\left( \mathcal{H}^{p}\Big|Z_{n}>0,Z_{0}=1\right) \Longrightarrow
\mathbb{L}^{+}\left( \mathcal{B}\right) ,
\end{equation*}%
where the symbol $\Longrightarrow $ stands for the weak convergence in the
space $D\left( [0,\infty ),\mathbb{R}\right) $ of c\`{a}dl\`{a}g functions
in $[0,\infty )$ endowed with the Skorokhod topology. In particular,
\begin{equation*}
\lim_{n\rightarrow \infty }\mathbf{P}\left( \frac{\log Z_{p}}{c_{p}}\leq z%
\Big|Z_{n}>0,Z_{0}=1\right) =\mathbb{P}^{+}\left( B_{1}\leq z\right) =C_{0}%
\mathbb{E}^{(m)}\left[ I\left\{ B_{1}\leq z\right\} B_{1}^{\alpha \left(
1-\rho \right) }\right]
\end{equation*}%
for any $z>0$.
\end{theorem}

\textbf{Remark 1}. This theorem complements Corollary 1.6 in \cite{AGKV05},
which states that under Conditions $A1$ and $A2$%
\begin{equation*}
\mathcal{L}\left( \mathcal{G}^{n}\Big|Z_{n}>0,Z_{0}=1\right) \Longrightarrow
\mathbb{L}^{(m)}\left( \mathcal{B}^{1}\right)
\end{equation*}%
as $n\rightarrow \infty ,$ where the symbol $\Longrightarrow $ stands for
the weak convergence in the space $D([0,1],\mathbb{R)}$ of c\`{a}dl\`{a}g
functions in $[0,1]$ endowed with the Skorokhod topology. In particular,
\begin{equation*}
\lim_{n\rightarrow \infty }\mathbf{P}\left( \frac{\log Z_{n}}{c_{n}}\leq z%
\Big|Z_{n}>0,Z_{0}=1\right) =\mathbb{P}^{(m)}\left( B_{1}\leq z\right) =%
\mathbb{E}^{(m)}\left[ I\left\{ B_{1}\leq z\right\} \right]
\end{equation*}%
for any $z>0$.

Let, for $U>0$
\begin{equation*}
\mathcal{H}_{U}^{p}:=\left\{ \frac{\log Z_{\left[ pu\right] }}{c_{p}},0\leq
u\leq U\right\} .
\end{equation*}

\begin{corollary}
\label{C_twoProc}If Condition $A$ is valid, then, for any $U>0$%
\begin{equation*}
\mathcal{L}\left( (\mathcal{H}_{U}^{p},\mathcal{G}^{n})\Big|%
Z_{n}>0,Z_{0}=1\right) \Longrightarrow \mathbb{L}^{+}\left( \mathcal{B}%
^{U}\right) \times \mathbb{L}^{(m)}\left( \mathcal{B}^{1}\right)
\end{equation*}
as $n\rightarrow \infty .$
\end{corollary}

\textbf{Remark 2}. If $\mathbf{E}X=0$ and $VarX\in \left( 0,\infty \right) $
then, for any $z>0$
\begin{equation*}
\mathbb{P}^{(m)}\left( B_{1}\leq z\right)
=\int_{0}^{z}xe^{-x^{2}/2}dx=1-e^{-z^{2}/2}\text{ }
\end{equation*}%
and%
\begin{equation*}
\text{ }\mathbb{P}^{+}\left( B_{1}\leq z\right) =\sqrt{\frac{2}{\pi }}%
\int_{0}^{z}x^{2}e^{-x^{2}/2}dx.
\end{equation*}

We have seen by (\ref{tails}) that the asymptotic behavior of the survival
probability of the process $\mathcal{Z}$ is primarily determined by the
random walk $\mathcal{S}$, since only the constant \ $\theta $ depends on
the fine structure of $\mathcal{Z}$ (see formula (\ref{DefTheta}) below).
However, one also has to take into account that the random walk changes its
properties drastically, when conditioned on the event $\{Z_{n}>0\}$. The
next theorem, describing the trajectories of the random walk $\mathcal{S}$
that provide survival of the critical process in a random environment at the
initial stage of the development of the population, illustrates this fact.

For $U\in (0,\infty ]$ let
\begin{eqnarray*}
\mathcal{Q}_{U}^{p} &:&=\left\{ \frac{S_{\left[ pu\right] }}{c_{p}},0\leq
u\leq U\right\} ,\quad \mathcal{Q}^{p}=\mathcal{Q}_{\infty }^{p}, \\
\mathcal{S}_{U}^{n} &:&=\left\{ \frac{S_{pU+\left[ \left( n-pU\right) t%
\right] }}{c_{n}},0\leq t\leq 1\right\} ,\quad \mathcal{S}^{n}:=\mathcal{S}%
_{0}^{n}.
\end{eqnarray*}

\begin{theorem}
\label{T_RWcondit}If Conditions $A$ is valid then, as $n\rightarrow \infty $
\begin{equation*}
\mathcal{L}\left( \mathcal{Q}^{p}\Big|Z_{n}>0,Z_{0}=1\right) \Longrightarrow
\mathbb{L}^{+}\left( \mathcal{B}\right) .
\end{equation*}
\end{theorem}

\textbf{Remark 3}. This theorem complements Theorem 1.5 in \cite{AGKV05},
which states that under Conditions $A1$ and $A2$%
\begin{equation}
\mathcal{L}\left( \mathcal{S}^{n}\Big|Z_{n}>0,Z_{0}=1\right) \Longrightarrow
\mathbb{L}^{(m)}\left( \mathcal{B}^{1}\right)   \label{FunctAGKV}
\end{equation}%
as $n\rightarrow \infty .$

\begin{corollary}
\label{C_TwoRW}If Condition $A$ is valid, then for any $U>0$
\begin{equation*}
\mathcal{L}\left( (\mathcal{Q}_{U}^{p},\mathcal{S}^{n})\Big|%
Z_{n}>0,Z_{0}=1\right) \Longrightarrow \mathbb{L}^{+}\left( \mathcal{B}%
^{U}\right) \times \mathbb{L}^{(m)}\left( \mathcal{B}^{1}\right)
\end{equation*}%
as $n\rightarrow \infty $.
\end{corollary}

The usage of the associated random walks to study branching processes in
random environment has a long history. It seems that Kozlov \cite{Koz76} was
the first who observed that to investigate properties of the critical
branching processes in random environment it is convenient to use ladder
epochs of the associated random walks. This fact has been used in various
situations for the case of the associated random walks with zero or negative
drift and finite variance of increments (see \cite{Af1993}, \cite{Af1997},%
\cite{Af1999a}, \cite{Af1999b}, \cite{Af2001} \cite{GK2000}, \cite{Koz1995}
and \cite{Vat2002}). The first steps to overcome the assumption of a finite
variance random walk in the driftless case were taken in \cite{DGV2004} and
\cite{VD2003}. In recent years papers \cite{AGKV05}, \cite{AGKV05b}, \cite%
{ABKV2009}, \cite{ABGK2013}, \cite{BDKV}, \cite{VL15} and some others
provide a systematic approach to the study of branching processes in random
environment under rather general assumptions on the properties of the
associated random walk (see, surveys \cite{Vat2016} and \cite{VDS2013} for a
detailed exposition).

\section{Auxiliary results}

We will use the symbols $K,K_{1},K_{2},...$ to denote different constants.
They are not necessarily the same in different formulas.

\subsection{Properties of the associated random walk}

To prove the main results of the pepar we need to know the asymptotic
behavior of the function $V(x)$ as $x\rightarrow \infty $. The following
lemma gives the desired asymptotics.

\begin{lemma}
\label{Renew2} (compare with Lemma 13 in \cite{VW09} and Corollary 8 in \cite%
{Don12}) If $X\in \mathcal{D}\left( \alpha ,\beta \right) $ then there
exists a slowly varying function $l_{0}(x)$ such that
\begin{equation}
V(x)\sim x^{\alpha (1-\rho )}l_{0}(x)  \label{CombRR1}
\end{equation}%
as $x\rightarrow \infty $.
\end{lemma}

Our next result is a combination (with a slight reformulation) of Lemma 2.1
in \cite{AGKV05} and Corollaries 3 and 8 in \cite{Don12}:

\begin{lemma}
\label{L_minimBelow} If $X\in \mathcal{D}\left( \alpha ,\beta \right) ,$
then there exist positive constants $K,K_{1}$ and $K_{2}$ such that, as\ $%
n\rightarrow \infty $
\begin{equation}
\mathbf{P}\left( L_{n}\geq -w\right) \sim V(w)\mathbf{P}\left( L_{n}\geq
0\right) \sim KV(w)n^{\rho -1}l(n)  \label{UnifPrecise}
\end{equation}%
uniformly for $0\leq w\ll c_{n},$ and%
\begin{equation}
\mathbf{P}\left( L_{n}\geq -w\right) \leq K_{1}V(w)n^{\rho -1}l(n)\leq
K_{2}V(w)\mathbf{P}\left( L_{n}\geq 0\right) ,\ w\geq 0,n\geq 1.
\label{IntermedEst}
\end{equation}
\end{lemma}

For further references we prove the following simple statement.

\begin{lemma}
\label{L_NewUniform}Let $\mathcal{A}_{n}\subset \mathbb{R},$ $n\in \mathbb{N}%
,$ be a family of subsets and let $b_{n}(x),$ $n\in \mathbb{N},$ be a
sequence of functions such that, for any fixed sequence $\left\{ a_{n},n\in
\mathbb{N}\right\} $ such that $a_{n}\in \mathcal{A}_{n}$ for all $n\in
\mathbb{N}$
\begin{equation}
\lim_{n\rightarrow \infty }b_{n}(a_{n})=0.  \label{Contr}
\end{equation}%
Then%
\begin{equation*}
\lim_{n\rightarrow \infty }\sup_{a\in \mathcal{A}_{n}}\left\vert
b_{n}(a)\right\vert =0.
\end{equation*}
\end{lemma}

\textbf{Proof}. Assume that the conclusion of the lemma is not true. Then,
there exists $\varepsilon >0$ such that for all $N$ there exist $n(N)\geq N$
and $a_{n(N)}\in \mathcal{A}_{n(N)}$ such that%
\begin{equation*}
\left\vert b_{n(N)}\left( a_{n(N)}\right) \right\vert \geq \varepsilon .
\end{equation*}

This, clearly, contradicts (\ref{Contr}).

The lemma is proved.

In the sequel we agree to consider the expressions of the form $\lim A(p,n)$
or $\limsup A(p,n)$ without lower indices as the $\lim $ or $\limsup $ of
the triangular array $\left\{ A(p,n),p\geq 1,n\geq 1\right\} $ calculated
under the assumption $pn^{-1}\rightarrow 0$ as $p,n\rightarrow \infty $. We
also write $a_{n}\ll b_{n}~$if $\lim_{n\rightarrow \infty }a_{n}/b_{n}=0$.

Let $\phi _{1}:\Omega _{1}\rightarrow \mathbb{R}$ be a bounded uniformly
continuous functional and $\left\{ \varepsilon _{n},n\in \mathbb{N}\right\} $
be a sequence of positive numbers vanishing as $n\rightarrow \infty $.

\begin{lemma}
\label{L_uniform}If Condition $A1$ is valid then
\begin{equation}
\mathbf{E}\left[ \phi _{1}\left( \mathcal{S}^{n}\right) |L_{n}\geq -x\right]
\rightarrow \mathbb{E}^{(m)}\left[ \phi _{1}(\mathcal{B}^{1})\right]
\label{DonDur}
\end{equation}%
as $n\rightarrow \infty $ uniformly in $0\leq x\leq \varepsilon _{n}c_{n}.$
\end{lemma}

\textbf{Proof of Lemma \ref{L_uniform}}. It was shown in Theorem 1.1 of
\cite{CCh2006} that, given Condition $A1$ convergence (\ref{DonDur}) holds
for any sequence $x=x_{n}$ meeting the restriction $0\leq x_{n}\ll c_{n}$ as
$n\rightarrow \infty $. This and Lemma \ref{L_NewUniform} \ with $\mathcal{A}%
_{n}:=\left\{ 0\leq x\leq \varepsilon _{n}c_{n}\right\} $ imply the desired
statement.

Now we are ready to demonstrate the validity of the following result.

\begin{lemma}
\label{L_RW}If Conditions $A1$ and (\ref{Msmalln}) are valid then, for $U>0$
and any $r\geq 0$
\begin{equation*}
\mathcal{L}\left( (\mathcal{Q}_{U}^{p},\mathcal{S}^{n})\Big|L_{n}\geq
-r\right) \Longrightarrow \mathbb{L}^{+}\left( \mathcal{B}^{U}\right) \times
\mathbb{L}^{(m)}\left( \mathcal{B}^{1}\right)
\end{equation*}%
as $n\rightarrow \infty $.
\end{lemma}

\textbf{Proof}. Consider the processes $\mathcal{S}^{k,n}$ and $\mathcal{%
\tilde{S}}^{k,n},0\leq k\leq n,$ given by%
\begin{equation}
S_{t}^{k,n}:=\frac{S_{\left[ nt\right] \wedge k}}{c_{n}},\quad \tilde{S}%
^{k,n}:=\frac{1}{c_{n}}\left( S_{\left[ nt\right] }-S_{\left[ nt\right]
\wedge k}\right) ,\text{ }0\leq t\leq 1.  \label{DefCens}
\end{equation}%
Clearly,%
\begin{equation*}
\mathcal{S}^{n}=\mathcal{S}^{k,n}+\mathcal{\tilde{S}}^{k,n}.
\end{equation*}%
Let $\mathcal{S}^{\ast }:=\left\{ S_{n}^{\ast },n\geq 0\right\} $ be a
probabilistic and independent copy of the random walk $\mathcal{S}=\left\{
S_{n},n\geq 0\right\} $ and
\begin{equation*}
L_{n}^{\ast }:=\min \left( S_{0}^{\ast },S_{1}^{\ast },...,S_{n}^{\ast
}\right) ,\quad \left( \mathcal{S}^{\ast }\right) _{U}^{n}:=\left\{ \frac{S_{%
\left[ \left( n-pU\right) t\right] }^{\ast }}{c_{n}},0\leq t\leq 1\right\} .
\end{equation*}%
For a fixed $N>0$ set%
\begin{equation*}
I_{N}(x):=\left\{
\begin{array}{ccc}
0 & \text{if} & x\leq N^{-1}, \\
Nx-1 & \text{if} & x\in \left( N^{-1},2N^{-1}\right) , \\
1 & \text{if} & 2N^{-1}\leq x\leq N, \\
N+1-x & \text{if} & N<x\leq N+1, \\
0 & \text{if} & x>N+1,%
\end{array}%
\right.
\end{equation*}%
and let
\begin{equation*}
\phi :\Omega _{U}\rightarrow \mathbb{R}\text{ and }\phi _{1}:\Omega
_{1}\rightarrow \mathbb{R}
\end{equation*}%
be two continuous and bounded functionals.

Then, for fixed positive $U$ and $N$ and $pU=n\varepsilon _{n}$, where $%
\varepsilon \geq \varepsilon _{n}\downarrow 0$ as $n\rightarrow \infty ,$\
we have (with a slight abuse of notation)%
\begin{eqnarray*}
&&\mathbf{E}\left[ \phi \left( \mathcal{Q}_{U}^{p}\right) I_{N}\left( \frac{%
S_{pU}}{c_{p}}\right) \phi _{1}\left( \mathcal{S}^{n}\right) ;L_{n}\geq -r%
\right]  \\
&=&\mathbf{E}\left[ \phi \left( \mathcal{Q}_{U}^{p}\right) I_{N}\left( \frac{%
S_{pU}}{c_{p}}\right) I\left\{ L_{pU}\geq -r\right\} \mathbf{E}\left[ \phi
_{1}\left( \left( \mathcal{S}^{\ast }\right) _{U}^{n}+\mathcal{S}%
^{pU,n}\right) I\left\{ L_{n-pU}^{\ast }\geq -S_{pU}-r\right\} \right] %
\right] .
\end{eqnarray*}%
Here and in what follows we agree to consider $pU$ and $n-pU$ as $\left[ pU%
\right] $ and $\left[ n-pU\right] ,$ respectively. Since $%
c_{p}/c_{n}\rightarrow 0$ as $n\rightarrow \infty ,$ it follows that, given $%
L_{pU}\geq -r$
\begin{equation*}
\frac{S_{pU}}{c_{n}}I_{N}\left( \frac{S_{pU}}{c_{p}}\right) \rightarrow 0%
\text{ a.s. }
\end{equation*}%
and $\mathcal{S}^{pU,n}$ vanishes as $n\rightarrow \infty $. This
observation, Lemma \ref{L_uniform} and the continuity of $\phi _{1}$ imply
\begin{equation*}
\mathbf{E}\left[ \phi _{1}\left( \left( \mathcal{S}^{\ast }\right) _{U}^{n}+%
\mathcal{S}^{pU,n}\right) |L_{n(1-\varepsilon _{n})}^{\ast }\geq -S_{pU}-r%
\right] \rightarrow \mathbb{E}^{(m)}\left[ \phi _{1}(\mathcal{B}^{1})\right]
\end{equation*}%
as $n\rightarrow \infty $ uniformly for $0\leq S_{pU}\leq Nc_{p}\ll c_{n}.$
On the other hand, by (\ref{UnifPrecise}), (\ref{AsH}), (\ref{CombRR1}) and
properties of regularly varying functions (see, for instance, \cite{Sen76})
we deduce, as $p,n\rightarrow \infty $:
\begin{eqnarray*}
\mathbf{P}\left( L_{n-pU}^{\ast }\geq -S_{pU}-r\right) I_{N}\left( \frac{%
S_{pU}}{c_{p}}\right)  &\sim &V\left( S_{pU}\right) I_{N}\left( \frac{S_{pU}%
}{c_{p}}\right) \mathbf{P}\left( L_{n}\geq 0\right)  \\
&=&\frac{V\left( S_{pU}\right) }{V(c_{p})}I_{N}\left( \frac{S_{pU}}{c_{p}}%
\right) \times V(c_{p})\mathbf{P}\left( L_{n}\geq 0\right)  \\
&\sim &\left( \frac{S_{pU}}{c_{p}}\right) ^{\alpha \left( 1-\rho \right)
}I_{N}\left( \frac{S_{pU}}{c_{p}}\right) \frac{C_{0}\mathbf{P}\left(
L_{n}\geq 0\right) }{\mathbf{P}\left( L_{p}\geq 0\right) } \\
&\sim &\left( \frac{S_{pU}}{c_{p}}\right) ^{\alpha \left( 1-\rho \right)
}I_{N}\left( \frac{S_{pU}}{c_{p}}\right) \frac{C_{0}\mathbf{P}\left(
L_{n}\geq -r\right) }{\mathbf{P}\left( L_{p}\geq -r\right) }.
\end{eqnarray*}%
Hence we get after evident but awkward transformations that, as $%
p,n\rightarrow \infty $%
\begin{eqnarray*}
&&\mathbf{E}\left[ \phi \left( \mathcal{Q}_{U}^{p}\right) \phi _{1}\left(
\mathcal{S}^{n}\right) I_{N}\left( \frac{S_{pU}}{c_{p}}\right) |L_{n}\geq -r%
\right]  \\
&&\qquad \sim C_{0}\mathbb{E}^{(m)}\left[ \phi _{1}(\mathcal{B}_{1})\right]
\mathbf{E}\left[ \phi \left( \mathcal{Q}_{U}^{p}\right) \left( \frac{S_{pU}}{%
c_{p}}\right) ^{\alpha \left( 1-\rho \right) }I_{N}\left( \frac{S_{pU}}{c_{p}%
}\right) |L_{pU}\geq -r\right] .
\end{eqnarray*}%
By Theorem 1.1 of \cite{CCh2006}, as $p\rightarrow \infty $
\begin{eqnarray*}
&&\mathbf{E}\left[ \phi \left( \mathcal{Q}_{U}^{p}\right) \left( \frac{S_{pU}%
}{c_{p}}\right) ^{\alpha \left( 1-\rho \right) }I_{N}\left( \frac{S_{pU}}{%
c_{p}}\right) |L_{pU}\geq -r\right]  \\
&&\qquad \qquad \qquad \rightarrow \mathbb{E}^{(m)}\left[ \phi \left(
\mathcal{B}^{U}\right) B_{U}^{\alpha \left( 1-\rho \right) }I_{N}\left(
B_{U}\right) \right] =\mathbb{E}^{+}\left[ \phi \left( \mathcal{B}%
^{U}\right) I_{N}\left( B_{U}\right) \right] .
\end{eqnarray*}%
Thus, under Conditions $A1$ and (\ref{Msmalln})%
\begin{eqnarray*}
&&\lim \mathbf{E}\left[ \phi \left( \mathcal{Q}_{U}^{p}\right) I_{N}\left(
\frac{S_{pU}}{c_{p}}\right) \phi _{1}\left( \mathcal{S}^{n}\right)
|L_{n}\geq -r\right]  \\
&&\qquad \qquad \qquad \qquad \qquad \qquad =C_{0}\mathbb{E}^{+}\left[ \phi
\left( \mathcal{B}^{U}\right) I_{N}\left( B_{U}\right) \right] \times
\mathbb{E}^{(m)}\left[ \phi _{1}(\mathcal{B}^{1})\right] .
\end{eqnarray*}%
Letting now $N\rightarrow \infty $ we get
\begin{equation*}
\mathcal{L}\left( (\mathcal{Q}_{U}^{p},\mathcal{S}^{n})\Big|L_{n}\geq
-r\right) \Longrightarrow \mathbb{L}^{+}\left( \mathcal{B}^{U}\right) \times
\mathbb{L}^{(m)}\left( \mathcal{B}^{1}\right)
\end{equation*}%
for any $U>0$. \textbf{\ }

The lemma is proved.

\begin{corollary}
\label{C_funct}If Conditions $A1$ and (\ref{Msmalln}) are valid then
\begin{equation*}
\mathcal{L}\left( \mathcal{Q}^{p}\Big|L_{n}\geq -r\right) \Longrightarrow
\mathbb{L}^{+}\left( \mathcal{B}\right)
\end{equation*}%
as $n\rightarrow \infty $.
\end{corollary}

\textbf{Proof}. It follows from Lemma \ref{L_RW} that
\begin{equation*}
\mathcal{L}\left( \mathcal{Q}_{U}^{p}\Big|L_{n}\geq -r\right)
\Longrightarrow \mathbb{L}^{+}\left( \mathcal{B}^{U}\right)
\end{equation*}%
for any $U>0$. This fact combined with Theorem 16.7 in \cite{Bil99}
completes the proof of the corollary.

\section{Conditional limit theorem}

For convenience we introduce the notation%
\begin{equation*}
A_{u.s.}=\left\{ Z_{n}>0\text{ for all }n\geq 0\right\}
\end{equation*}%
and recall that by Corollary 1.2 in \cite{AGKV05}, (\ref{tails}) and~(\ref%
{AsH})
\begin{equation}
\mathbf{P}\left( Z_{n}>0\right) \sim \theta \mathbf{P}\left( L_{n}\geq
0\right) \sim \theta n^{-(1-\rho )}l(n)\sim \frac{\theta C_{0}}{V(c_{n})}
\label{DonRen}
\end{equation}%
as $n\rightarrow \infty $, where
\begin{equation}
\theta =\sum_{k=0}^{\infty }\mathbf{E}[\mathbf{P}_{Z_{k}}^{+}\left(
A_{u.s.}\right) ;\tau _{k}=k].  \label{DefTheta}
\end{equation}%
Let%
\begin{equation*}
\hat{L}_{k,n}:=\min_{0\leq j\leq n-k}\left( S_{k+j}-S_{k}\right)
\end{equation*}%
and let $\mathcal{\tilde{F}}_{k}$ be the $\sigma -$algebra generated by the
tuple $\left\{ Z_{0},Z_{1},...,Z_{k};Q_{1},Q_{2},...,Q_{k}\right\} $ (see (%
\ref{DefEnvir})). For further references we formulate two statements borrowed from \cite{AGKV05}.

\begin{lemma}
\label{L_cond}(see Lemma 2.5 in \cite{AGKV05}) Assume Condition $A1$. Let $%
Y_{1},Y_{2},...$be a uniformly bounded sequence of real-valued random
variables adapted to the filtration $\mathcal{\tilde{F}=}\left\{ \mathcal{%
\tilde{F}}_{k},k\in \mathbb{N}\right\} $, which converges $\mathbf{P}^{+}$%
-a.s. to some random variable $Y_{\infty }$. Then, as $n\rightarrow \infty $%
\begin{equation*}
\mathbf{E}\left[ Y_{n}|L_{n}\geq 0\right] \rightarrow \mathbf{E}^{+}\left[
Y_{\infty }\right] .
\end{equation*}
\end{lemma}

Denote%
\begin{equation}
\tau _{n}:=\min \left\{ j:S_{j}=L_{n}\right\} .  \label{Deftau}
\end{equation}

\begin{lemma}
\label{generaltheorem}(see Lemma 4.1 in \cite{AGKV05}) Assume Conditions A1
and let $l\in \mathbb{N}_{0}$. Suppose that $\zeta _{1},\zeta _{2},...$ is a
uniformly bounded sequence of real-valued random variables, which, for every
$k\geq 0$ meets the equality%
\begin{equation}
\mathbf{E}\left[ \zeta _{n};Z_{k+l}>0,\hat{L}_{k,n}\geq 0\ |\mathcal{\tilde{F%
}}_{k}\right] =\mathbf{P}\left( L_{n}\geq 0\right) \left( \zeta _{k,\infty
}+o(1)\right) ,\qquad \mathbf{P}\text{-a.s.}  \label{Start}
\end{equation}%
as $n\rightarrow \infty $ with random variables $\zeta _{1,\infty }=\zeta
_{1,\infty }\left( l\right) ,\zeta _{k,\infty }=\zeta _{2,\infty }\left(
l\right) ,....$ Then%
\begin{equation*}
\mathbf{E}\left[ \zeta _{n};Z_{\tau _{n}+l}>0\right] =\mathbf{P}\left(
L_{n}\geq 0\right) \left( \sum_{k=0}^{\infty }\mathbf{E}\left[ \zeta
_{k,\infty };\tau _{k}=k\right] +o(1)\right)
\end{equation*}%
as $n\rightarrow \infty $, where the right-hand side series is absolutely
convergent.
\end{lemma}

For $U>0$ and $q\leq p,\quad pU\leq n$ let%
\begin{eqnarray*}
\mathcal{X}_{U}^{q,p} &:&=\left\{ X_{u}^{q,p}=e^{-S_{q+\left[ u(p-q)\right]
}}Z_{q+\left[ u(p-q)\right] },0\leq u\leq U\right\} ,\quad \\
\mathcal{X}^{q,p} &:&=\left\{ X_{u}^{q,p}=e^{-S_{q+\left[ u(p-q)\right]
}}Z_{q+\left[ u(p-q)\right] },0\leq u<\infty \right\} , \\
\mathcal{Y}_{U}^{p,n} &:&=\left\{ Y_{t}^{p,n}=e^{-S_{pU+\left[ (n-pU)t\right]
}}Z_{pU+\left[ (n-pU)t\right] },0\leq t\leq 1\right\} ,\quad \mathcal{Y}%
^{p,n}:=\mathcal{Y}_{0}^{p,n}.
\end{eqnarray*}

The next statement is an evident corollary of Theorem 1.3 in \cite{AGKV05}
and we give its proof for completeness only.

\begin{lemma}
\label{L_skoroh}Assume Conditions $A1$ and $A2$. Let $\left(
q_{1},p_{1}\right) ,\left( q_{2},p_{2}\right) ,...$ be a sequence of pairs
of positive integers such that $q_{n}\ll p_{n}$ as $n\rightarrow
\infty $. If $p_{n}$ $\ll n$ then, for any $U>0$%
\begin{equation*}
\mathcal{L}\left( \left( \mathcal{X}_{U}^{q_{n},p_{n}},\mathcal{Y}%
_{U}^{p_{n},n}\right) \ |Z_{n}>0,Z_{0}=1\right) \Longrightarrow \mathcal{L}%
\left( (W_{u},\ 0\leq u\leq U),(\breve{W}_{t},\ 0\leq t\leq 1)\right) \
\end{equation*}%
as $n\rightarrow \infty $, where
\begin{equation}
\mathbf{P}\left( W_{u}=\breve{W}_{t}=W,\ 0\leq u\leq U,0\leq t\leq 1\right)
=1
\end{equation}%
for some random variable $W$ such that%
\begin{equation*}
\mathbf{P}\left( 0<W<\infty \right) =1.
\end{equation*}
\end{lemma}

\textbf{Proof}. We follow (with minor changes) the line of proving Theorem
1.3 in~\cite{AGKV05}. According to Proposition 3.1 in \cite{AGKV05} there
exists a strictly positive and finite random variable \ $W^{+}$ such that,
as $n\rightarrow \infty $%
\begin{equation}
e^{-S_{n}}Z_{n}\rightarrow W^{+}\text{ }\quad \mathbf{P}^{+}\text{--}a.s.
\label{Prop31}
\end{equation}%
and%
\begin{equation}
\left\{ W^{+}>0\right\} =\left\{ Z_{n}>0\text{ for all }n\right\} \text{ }%
\quad \mathbf{P}^{+}\text{--}a.s.  \label{Prop31b}
\end{equation}%
Fix $U>0$ and let $\phi $ be a bounded continuous function on the space $%
\Omega _{U}=D\left( [0,U],\mathbb{R}\right) $ of c\`{a}dl\`{a}g functions
and let $\phi _{1}$ be a bounded continuous function on the space $\Omega
_{1}$. For $s\in \mathbb{R}$ let $\mathcal{W}_{U}^{s}:=\left\{
W_{u}^{s},0\leq u\leq U\right\} $ and $\mathcal{\check{W}}^{s}:=\left\{
\check{W}_{t}^{s},0\leq t\leq 1\right\} $ denote the processes with constant
paths coinciding (formally) within the time-interval $\left[ 0,\min \left\{
U,1\right\} \right] $:
\begin{equation*}
W_{u}^{s}\ :=\ e^{-s}W^{+}\ ,\quad 0\leq u\leq U,\quad \check{W}_{t}^{s}:=\
e^{-s}W^{+},\quad 0\leq t\leq 1.
\end{equation*}%
It follows from (\ref{Prop31}) that, for fixed $s\in \mathbb{R}$ the
two-dimensional process
$$(e^{-s}\mathcal{X}_{U}^{q_{n},p_{n}},e^{-s}\mathcal{%
Y}_{U}^{p_{n},n})$$
converges, as $n,p_{n}\rightarrow \infty $ with $%
q_{n}\leq p_{n}\ll n$, to $\left( \mathcal{W}_{U}^{s},\mathcal{\check{W}}%
^{s}\right) $ in the metric of uniform convergence and, consequently, in the
Skorokhod metric on the space $\Omega _{U}\times \Omega _{1}$ $\mathbf{P}^{+}
$--a.s.$,$ and%
\begin{eqnarray*}
\mathcal{K}_{n}&:=&\phi (e^{-s}\mathcal{X}_{U}^{q_{n},p_{n}})\phi
_{1}(e^{-s}\mathcal{Y}_{U}^{p_{n},n})I\left( Z_{n}>0\right) \  \\
&\rightarrow & \mathcal{K}_{\infty }:=\phi (\mathcal{W}_{U}^{s})\phi _{1}(%
\mathcal{\check{W}}^{s})I\{W^{+}>0\}\text{ \ \ \ }\mathbf{P}^{+}-\text{a.s.}
\end{eqnarray*}%
For $q\leq p\leq n$ and $z\in \mathbb{N}_{0}$ define
\begin{eqnarray*}
\psi (z,s,q,p,n)&:=&\mathbf{E}_{z}[\phi (e^{-s}\mathcal{X}_{U}^{q,p})\phi
_{1}(e^{-s}\mathcal{Y}_{U}^{p,n});Z_{n}>0,L_{n}\geq 0]\  \\
&=&\mathbf{E}_{z}[\phi (e^{-s}\mathcal{X}_{U}^{q,p})\phi _{1}(e^{-s}\mathcal{%
Y}_{U}^{p,n})I\left( Z_{n}>0\right) |L_{n}\geq 0]\mathbf{P}\left( L_{n}\geq
0\right) .
\end{eqnarray*}%
Since $\mathcal{K}_{n}\ \rightarrow \ \mathcal{K}_{\infty }$\ \ $\mathbf{P}%
^{+}-$a.s. as $n\rightarrow \infty $, it follows from Lemma \ref{L_cond}
that
\begin{equation*}
\psi (z,s,q_{n},p_{n},n)\ =\mathbf{P}\left( L_{n}\geq 0\right) \;\big(%
\mathbf{E}_{z}^{+}[\phi (\mathcal{W}_{U}^{s})\phi _{1}(\mathcal{\check{W}}%
^{s});W^{+}>0]+o(1)\big)\ .
\end{equation*}%
Observe now that, for $k\leq q\leq p\leq n$
\begin{equation*}
\mathbf{E}[\phi (e^{-s}\mathcal{X}_{U}^{q,p})\phi _{1}(e^{-s}\mathcal{Y}%
_{U}^{p,n});Z_{n}>0,\hat{L}_{k,n}\geq 0\;|\;\mathcal{F}_{k}]\ =\ \psi
(Z_{k},S_{k},q-k,p-k,n-k).
\end{equation*}%
Therefore, we may apply Lemma \ref{generaltheorem} to the
random variables
\begin{equation*}
\zeta _{n}=\phi (e^{-s}\mathcal{X}_{U}^{q_{n},p_{n}})\phi _{1}(e^{-s}%
\mathcal{Y}_{U}^{p_{n},n})I\{Z_{n}>0\}
\end{equation*}%
and
\begin{equation*}
\zeta _{k,\infty }=\mathbf{E}_{Z_{k}}^{+}[\phi (\mathcal{W}_{U}^{S_{k}})\phi
_{1}(\mathcal{\check{W}}^{S_{k}});W^{+}>0]
\end{equation*}%
with $l=0$.

Using (\ref{DonRen}) we get
\begin{equation*}
\mathbf{E}[\phi (\mathcal{X}_{U}^{q_{n},p_{n}})\phi _{1}(\mathcal{Y}%
_{U}^{p_{n},n})\;|\;Z_{n}>0]\ \rightarrow \ \int \phi (\mathfrak{m})\phi
_{1}(\mathfrak{n})\,\lambda (d\mathfrak{m\times }d\mathfrak{n})\ \text{\ as
\ }n\rightarrow \infty ,
\end{equation*}%
where $\lambda $ is the measure on the product space of c\`{a}dl\`{a}g
functions on $\Omega _{U}\times \Omega _{1}$ specified by
\begin{equation*}
\lambda (d\mathfrak{m\times }d\mathfrak{n})\ :=\ \frac{1}{\theta }%
\sum_{k=0}^{\infty }\mathbf{E}[\lambda _{Z_{k},S_{k}}(d\mathfrak{m\times }d%
\mathfrak{n});Z_{k}>0,\tau _{k}=k]
\end{equation*}%
with
\begin{equation*}
\lambda _{z,s}(d\mathfrak{m\times }d\mathfrak{n})\ :=\ \mathbf{P}_{z}^{+}[%
\mathcal{W}_{U}^{s}\in d\mathfrak{m},\mathcal{\check{W}}^{s}\in d\mathfrak{n,%
}W^{+}>0].
\end{equation*}%
By (\ref{Prop31b}) the total mass of $\lambda _{z,s}$ is equal to $\mathbf{P}%
_{z}^{+}(Z_{n}>0$ for all $n\geq 0)$. Therefore, the representation of $%
\theta $ in (\ref{DefTheta}) shows that $\lambda $ is a probability measure.
Again using~(\ref{Prop31b}) we see that $\lambda _{z,s}$ is concentrated on
strictly positive constant functions only. Hence, the same is true for the
measure $\lambda $.

Lemma \ref{L_skoroh} is proved.

\begin{corollary}
\label{C_onedim} Assume Conditions $A1$ and $A2$. Let $\left(
q_{1},p_{1}\right) ,\left( q_{2},p_{2}\right) ,...$ be a sequence of pairs
of positive integers such that $q_{n}\ll p_{n}\ll n$ and $q_{n}\rightarrow
\infty $ as $n\rightarrow \infty $. Then
\begin{equation*}
\mathcal{L}\left( \mathcal{X}^{q_{n},p_{n}}\ |Z_{n}>0,Z_{0}=1\right)
\Longrightarrow \mathcal{L}\left( \left\{ W_{u},\ 0\leq u<\infty \right\}
\right) .
\end{equation*}
\end{corollary}

\textbf{Proof.} We know that%
\begin{equation*}
\mathcal{L}\left( \mathcal{X}_{U}^{q_{n},p_{n}}\ |Z_{n}>0,Z_{0}=1\right)
\Longrightarrow \mathcal{L}\left( \left\{ W_{u},\ 0\leq u\leq U\right\}
\right) \
\end{equation*}%
as $n\rightarrow \infty. $ for any $U>0$. This and Theorem 16.7 of \cite%
{Bil99} complete the proof of the corollary.

\textbf{Proof} \textbf{of Theorem \ref{T_RWcondit}}. Let $U>0$ be fixed. Consider the
processes
\begin{equation*}
\mathcal{Q}_{U}^{q,p}=\left\{ S_{u}^{q,p},0\leq u\leq U\right\} ,\mathcal{%
\tilde{Q}}_{U}^{q,p}=\left\{ \tilde{S}_{u}^{q,p},0\leq u\leq U\right\}
,\quad 0\leq q\leq pU,
\end{equation*}%
given by%
\begin{equation}
S_{u}^{q,p}:=\frac{S_{\left[ pu\right] \wedge q}}{c_{p}},\quad \tilde{S}%
_{u}^{q,p}:=\frac{1}{c_{p}}\left( S_{\left[ pu\right] }-S_{\left[ pu\right]
\wedge q}\right) ,\text{ }0\leq u\leq U.  \label{DefProc}
\end{equation}%
Clearly,
\begin{equation*}
\mathcal{Q}_{U}^{p}:=\mathcal{Q}_{U}^{q,p}+\mathcal{\tilde{Q}}_{U}^{q,p}.
\end{equation*}

Take $k,l\geq 0$ with $k+l\leq pU$. We may decompose the stochastic process $%
\mathcal{Q}_{U}^{p}$ as%
\begin{equation*}
\mathcal{Q}_{U}^{p}:=\mathcal{Q}_{U}^{k+l,p}+\mathcal{\tilde{Q}}_{U}^{k+l,p}.
\end{equation*}%
Let $\phi $ be a bounded continuous functional on $\Omega _{U}.$ Define
\begin{equation*}
\psi (\mathfrak{m},r):=\mathbf{E}[\phi (\mathfrak{m}+\mathcal{\tilde{Q}}%
_{U}^{k+l,p});\hat{L}_{k+l,n}\geq -r]
\end{equation*}%
for $\mathfrak{m}\in D[0,U]$ and $r\geq 0$. If $p,n\rightarrow \infty $ in
such a way that $pn^{-1}\rightarrow 0$ then, according to Corollary \ref%
{C_funct}
\begin{equation*}
\mathcal{L}\left( \left\{ S_{u}^{k+l,p},0\leq u<\infty \right\} \Big|\hat{L}%
_{k+p,n}\geq -r\right) \Longrightarrow \mathbb{L}^{+}\left( \left\{
B_{u},0\leq u<\infty \right\} \right)
\end{equation*}%
for each fixed pair $k$ and $l$. Hence, if the c\`{a}dl\`{a}g functions $%
\mathfrak{m}^{p}\in \Omega _{U}$ converge uniformly to the zero function as $%
p\rightarrow \infty $, then, given (\ref{Msmalln})
\begin{align*}
\psi (\mathfrak{m}^{p},r)\ & =\ \mathbf{P}\left( L_{n-(k+l)}\geq -r\right) (%
\mathbb{E}\,^{+}\left[ \phi (\mathcal{B}^{U})\right] +o(1)) \\
& =\ V(r)\mathbf{P}\left( L_{n}\geq 0\right) (\mathbb{E}\,^{+}\left[ \phi (%
\mathcal{B}^{U})\right] +o(1)),
\end{align*}%
as $p,n\rightarrow \infty $, where for the second equality we have applied (%
\ref{UnifPrecise}). Using the representation
\begin{equation}
\{\hat{L}_{k,n}\geq 0\}=\{\hat{L}_{k,k+l}\geq 0\}\cap \{\hat{L}_{k+l,n}\geq
-(S_{k+l}-S_{k})\}  \label{deco}
\end{equation}%
and taking into account that $\mathcal{Q}_{U}^{k+l,p}$ converges uniformly to zero $\mathbf{P}$%
--a.s. as $p\to\infty $, we have under Condition A:
\begin{eqnarray}
&&\mathbf{E}\left[ \,\phi \left( \mathcal{Q}_{U}^{p}\right) ;Z_{k+l}>0,\hat{L%
}_{k,n}\geq 0\;|\;\mathcal{F}_{k+l}\right]   \notag \\
&=&\psi \left( \mathcal{Q}_{U}^{k+l,p},S_{k+l}-S_{k}\right) I\left\{
Z_{k+l}>0,\,\hat{L}_{k,k+l}\geq 0\right\}   \notag \\
&=&V(S_{k+l}-S_{k})\mathbf{P}\left( L_{n}\geq 0\right) (\mathbb{E}\,^{+}%
\left[ \phi (\mathcal{B}^{U})\right] +o(1))I\left\{ Z_{k+l}>0,\,\hat{L}%
_{k,k+l}\geq 0\right\} \quad \mathbf{P}\text{--a.s}.  \notag \\
&&  \label{oct}
\end{eqnarray}%
This representation combined with (\ref{IntermedEst}) and (\ref{deco})
allows us to deduce the chain of estimates
\begin{eqnarray*}
&&\left\vert \mathbf{E}[\phi (\mathcal{Q}_{U}^{p});\ Z_{k+l}>0,\hat{L}%
_{k,n}\geq 0\;|\;\mathcal{F}_{k+l}]\right\vert \ \leq \ \sup |\phi |\,%
\mathbf{P}\left( \hat{L}_{k,n}\geq 0\;|\;\mathcal{F}_{k+l}\right)  \\
&&\qquad \qquad \qquad =\ \sup |\phi |\,\mathbf{P}\left( \hat{L}_{k+l,n}\geq
-(S_{k+l}-S_{k})\;|\;\mathcal{F}_{k+l}\right) I\left\{ \hat{L}_{k,k+l}\geq
0\right\}  \\
&&\qquad \qquad \qquad \leq \ K_{1}V(S_{k+l}-S_{k})\mathbf{P}\left(
L_{n-(k+l)}\geq 0\right) I\left\{ \hat{L}_{k,k+l}\geq 0\right\} \quad
\mathbf{P}\text{--a.s}.\
\end{eqnarray*}%
for some $K_{1}>0$. Observe now that according to (\ref{harm})
\begin{equation*}
\mathbf{E}[V(S_{k+l}-S_{k});\hat{L}_{k,k+l}\geq 0\;|\;\mathcal{F}%
_{k}]=V(0)<\infty \text{ \ \ \ \ \ }\mathbf{P}-\text{a.s.}
\end{equation*}%
Hence, using the dominated convergence theorem, (\ref{MesPm}) and the
definition of $\mathbf{P}^{+}$, we obtain by (\ref{oct}) that
\begin{align*}
\mathbf{E}[\phi \left( \mathcal{Q}_{U}^{p}\right) ;Z_{k+l}& >0,\hat{L}%
_{k,n}\geq 0\;|\;\mathcal{F}_{k}]\ =\ (\mathbb{E}\,^{+}\left[ \phi (\mathcal{%
B}^{U})\right] +o(1))\mathbf{P}\left( L_{n}\geq 0\right)  \\
& \qquad \qquad \times \mathbf{E}[V(S_{k+l}-S_{k});Z_{k+l}>0,\hat{L}%
_{k,k+l}\geq 0\;|\;\mathcal{F}_{k}] \\
& =\ (\mathbb{E}\,^{+}\left[ \phi (\mathcal{B}^{U})\right] +o(1))\mathbf{P}%
\left( L_{n}\geq 0\right) \mathbf{P}_{Z_{k}}^{+}\left( Z_{l}>0\right) \quad
\mathbf{P}\text{--a.s}.
\end{align*}%
Applying Lemma \ref{generaltheorem} to $\zeta _{n}=\phi (\mathcal{Q}%
_{U}^{p})$ with $n\gg p=p(n)\rightarrow \infty $ yields
\begin{eqnarray*}
&&\mathbf{E}[\phi \left( \mathcal{Q}_{U}^{p}\right) ;\ Z_{\tau _{n}+l}>0] \\
&&\qquad \qquad =(\mathbb{E}\,^{+}\left[ \phi (\mathcal{B}^{U})\right] +o(1))%
\mathbf{P}\left( L_{n}\geq 0\right) \sum_{k=0}^{\infty }\mathbf{E}[\mathbf{P}%
_{Z_{k}}^{+}\left( Z_{l}>0\right) ;\tau _{k}=k].
\end{eqnarray*}%
Therefore,
\begin{equation}
\mathbf{P}\left( Z_{\tau _{n}+l}>0\right) \ \sim \ \mathbf{P}\left(
L_{n}\geq 0\right) \sum_{k=0}^{\infty }\mathbf{E}[\mathbf{P}%
_{Z_{k}}^{+}\left( Z_{l}>0\right) ;\tau _{k}=k]  \label{wow2}
\end{equation}%
as $n\rightarrow \infty $, where the right-hand side series is convergent.
Observe that
\begin{eqnarray*}
&&\left\vert \mathbb{E}\,^{+}\left[ \phi (\mathcal{B}^{U})\right] \mathbf{P}%
\left( Z_{n}>0\right) -\mathbf{E}[\phi \left( \mathcal{Q}_{U}^{p}\right)
;Z_{n}>0]\right\vert  \\
&&\quad \leq \left\vert \mathbb{E}\,^{+}\left[ \phi (\mathcal{B}^{U})\right]
\mathbf{P}\left( Z_{n}>0\right) -\mathbf{E}[\phi \left( \mathcal{Q}%
_{U}^{p}\right) ;Z_{\tau _{n}+l}>0]\right\vert \hspace*{2cm} \\[-0.5ex]
&&\hspace*{2.7cm}\;+\ \sup |\phi |\;\mathbf{E}|I\left\{ Z_{n}>0\right\}
-I\left\{ Z_{{\tau _{n}}+l}>0\right\} |
\end{eqnarray*}%
and
\begin{eqnarray*}
\mathbf{E}|I\left\{ Z_{n}>0\right\} -I\left\{ Z_{\tau _{n}+l}>0\right\} |
&\leq &(\mathbf{P}\left( Z_{n}>0\right) -\mathbf{P}\left( Z_{n+l}>0\right) )
\\
&&+(\mathbf{P}\left( Z_{\tau _{n}+l}>0\right) -\mathbf{P}\left(
Z_{n+l}>0\right) ).
\end{eqnarray*}%
These estimates and (\ref{DonRen}) lead to the inequality
\begin{eqnarray}
&&\left\vert \mathbb{E}\,^{+}\left[ \phi (\mathcal{B}^{U})\right] -\mathbf{E}%
[\phi \left( \mathcal{Q}_{U}^{p}\right) \,|\,Z_{n}>0]\right\vert   \notag \\
&&\qquad \leq 2\sup |\phi |\left( \frac{1}{\theta }\sum_{k=0}^{\infty }%
\mathbf{E}[\mathbf{P}_{Z_{k}}^{+}\left( Z_{l}>0\right) ;\tau
_{k}=k]\;-1\right) \,+\,\varepsilon \left( p,n\right) ,  \label{wow}
\end{eqnarray}%
where $\lim \varepsilon \left( p,n\right) =0$. By the dominated convergence
theorem and the definition of \ $\theta $ in (\ref{DefTheta}) we conclude
that
\begin{equation*}
\sum_{k=0}^{\infty }\mathbf{E}[\mathbf{P}_{Z_{k}}^{+}\left( Z_{l}>0\right)
;\tau _{k}=k]\ \downarrow \ \theta \ \text{\ as }\,l\rightarrow \infty .
\end{equation*}%
Since the left-hand side of (\ref{wow}) does not depend on $l$, this gives
the assertion of Theorem \ref{T_RWcondit} for an arbitrary interval $0\leq
u\leq U$. To complete the proof of the theorem it remains to apply Theorem
16.7 of \cite{Bil99}.

\textbf{Proof of Corollary }\ref{C_TwoRW}. We use the notation of Lemma \ref%
{L_RW} and define%
\begin{equation*}
\psi ^{\ast }(\mathfrak{m,n,}r)\ :=\ \mathbf{E}[\phi (\mathfrak{m}+\mathcal{%
\tilde{Q}}_{U}^{k+l,p})\phi _{1}(\mathfrak{n}+\mathcal{\tilde{S}}^{k+l,p});%
\hat{L}_{k+l,n}\geq -r]
\end{equation*}%
for $\left( \mathfrak{m,n}\right) \in \Omega _{U}\times \Omega _{1}$ and $%
r\geq 0$. If a two dimensional vector of c\`{a}dl\`{a}g functions $\left(
\mathfrak{m}^{p}\mathfrak{,n}^{n}\right) \in \Omega _{U}\times \Omega _{1}$
converges uniformly to the two dimensional vector of zero
functions as $p=p(n)\rightarrow \infty $ as $n\rightarrow \infty ,$ and
condition (\ref{Msmalln}) is valid then, according to Lemma \ref{L_RW}%
\begin{eqnarray*}
\psi ^{\ast }(\mathfrak{m}^{p}\mathfrak{,n}^{n}\mathfrak{,}r)\  &=&\mathbf{P}%
\left( L_{n-(k+l)}\geq -r\right) (\mathbb{E}\,^{+}\left[ \phi (\mathcal{B}%
^{U})\right] \times \mathbb{E}\,^{(m)}\left[ \phi _{1}(\mathcal{B}^{1})%
\right] +o(1)) \\
&=&V(r)\mathbf{P}\left( L_{n}\geq 0\right) (\mathbb{E}\,^{+}\left[ \phi (%
\mathcal{B}^{U})\right] \times \mathbb{E}\,^{(m)}\left[ \phi _{1}(\mathcal{B}%
^{1})\right] +o(1)).
\end{eqnarray*}%
Let $k$ and $l$  be fixed. We know that the pair $\left( \mathcal{Q}%
^{k+l,p},S^{k+l,n}\right) $ uniformly converges, as $p,n\rightarrow \infty $ to the two dimensional
vector of zero functions $\mathbf{P}$--a.s. Hence we
obtain
\begin{eqnarray*}
&&\mathbf{E}\left[ \,\phi \left( \mathcal{Q}_{U}^{p}\right) \phi _{1}\left(
S^{n}\right) ;Z_{k+l}>0,\hat{L}_{k,n}\geq 0\;|\;\mathcal{F}_{k+l}\right]  \\
&&\qquad =\psi ^{\ast }\left( \mathcal{Q}_{U}^{k+l,p},\mathcal{S}%
^{k+l,n},S_{k+l}-S_{k}\right) I\left\{ Z_{k+l}>0,\,\hat{L}_{k,k+l}\geq
0\right\}  \\
&&\qquad =V(S_{k+l}-S_{k})\mathbf{P}\left( L_{n}\geq 0\right) \times (%
\mathbb{E}\,^{+}\left[ \phi (\mathcal{B}^{U})\right] \times \mathbb{E}%
\,^{(m)}\left[ \phi _{1}(\mathcal{B}^{1})\right] +o(1)) \\
&&\qquad \qquad \qquad \qquad \qquad \qquad \qquad \times I\left\{
Z_{k+l}>0,\,\hat{L}_{k,k+l}\geq 0\right\} \quad \mathbf{P}\text{--a.s}.
\end{eqnarray*}%
Repeating now almost literally (with evident changes) the proof of Theorem %
\ref{T_RWcondit} one can check the validity of Corollary \ref{C_TwoRW}.

\textbf{Proof of Theorem \ref{T_small}}. For each $U>0$ we have
\begin{eqnarray*}
&&\mathcal{L}\left( \left\{ \frac{\log Z_{q+up}}{c_{p}},0\leq u\leq
U\right\} \ \Big|Z_{n}>0,Z_{0}=1\right) \\
&&\qquad\quad =\mathcal{L}\left( \left\{ \frac{\log X_{u}^{q,p}}{c_{p}}+%
\frac{S_{pu}}{c_{p}},0\leq u\leq U\ \right\} \Big|Z_{n}>0,Z_{0}=1\right) .
\end{eqnarray*}%
This equality, Theorem~\ref{T_RWcondit} and Lemma \ref{L_skoroh} combined
with Theorem 16.7 of \cite{Bil99} justify the desired statement.

\textbf{Proof of Corollary }\ref{C_twoProc}. The needed statement follows
from the representation%
\begin{eqnarray*}
&&\mathcal{L}\left( \left\{ \frac{\log Z_{q+up}}{c_{p}},0\leq u\leq U;\frac{%
\log Z_{pU+\left[ \left( n-pU\right) t\right] }}{c_{n}},0\leq t\leq
1\right\} \ \Big|Z_{n}>0,Z_{0}=1\right) \\
&&=\mathcal{L}\left( \left\{ \frac{S_{pu}+\log X_{u}^{q,p}}{c_{p}},0\leq
u\leq U;\frac{S_{pU+\left[ \left( n-pU\right) t\right] }+\log Y_{t}^{p,n}}{%
c_{n}},0\leq t\leq 1\right\} \ \Big|Z_{n}>0,Z_{0}=1\right) ,
\end{eqnarray*}%
Lemma \ref{L_skoroh} and Corollary \ref{C_TwoRW}.


\begin{thebibliography}{99}
\bibitem{Af1993} Afanasyev V.I. A limit theorem for a critical branching
process in random environment. - \textit{Discrete Math. Appl.}, \textbf{5},
(1993), 45--58. (In Russian.)

\bibitem{Af1997} Afanasyev V.I. A new theorem for a critical branching
process in random environment. - \textit{Discrete Math. Appl.}, \textbf{7}
(1997), 497--513.

\bibitem{Af1999a} Afanasyev V.I. On the time of reaching a fixed level by a
critical branching process in a random environment. - \textit{Discrete Math.
Appl.}, \textbf{9} (1999), 627--643

\bibitem{Af1999b} Afanasyev V.I. On the maximum of a critical branching
process in a random environment. - \textit{Discrete Math. Appl.}, \textbf{9}
(1999), 267--284.

\bibitem{Af2001} Afanasyev V.I. A functional limit theorem for a critical
branching process in a random environment. - \textit{Discrete Math. Appl.},
\textbf{11} (2001), 587--606.

\bibitem{AGKV05} Afanasyev V.I., Geiger J., Kersting G., Vatutin V.A.
Criticality for branching processes in random environment. - \textit{Ann.
Probab.}, \textbf{33} (2005), 645--673.

\bibitem{AGKV05b} Afanasyev V.I., Geiger J., Kersting G., Vatutin V.A.
Functional limit theorems for strongly subcritical branching processes in
random environment. - \textit{Stoch. Proc. Appl.,} \textbf{115} (2005),
1658--1676.

\bibitem{ABKV2009} Afanasyev V.I., Boeinghoff Ch., Kersting G., Vatutin V.A.
Limit theorems for weakly subcritical branching processes in random
environment. - \textit{J. Theoret. Probab.}, \textbf{25 }(2012), 703--732.

\bibitem{ABGK2013} Afanasyev V.I., Boeinghoff Ch., Kersting G., Vatutin
V.A. Conditional limit theorems for intermediately subcritical branching
processes in random environment. - \textit{Ann. Inst. Henri-Poincar\'{e}},
\textbf{50} (2014), 602--627.

\bibitem{Bil99} Billingsley P. \textit{Convergence of Probability Measures}.
Willey, New York-London-Sydney-Toronto, 2nd ed., 1999.


\bibitem{BDKV} Beoinghoff C., Dyakonova E.E., Kersting G., and Vatutin V.A.
Branching processes in random environment which extinct at a given moment.
-- \textit{Markov Process. Relat. Fields}, \textbf{16} (2010), 329-350.

\bibitem{Chau97} Chaumont L. Conditionings and path decompositions for Levy
processes. - \textit{Stochastic Process. Appl.,} \textbf{64} (1996), 39--54.

\bibitem{Chau06} Chaumont L. Excursion normalisee, meandre at pont pour les
processus de Levy stables. - \textit{Bull. Sci. Math.}, \textbf{121} (1997),
5, 377--403.

\bibitem{CCh2006} Caravenna F., Chaumont L. Invariance principles for
random walks conditioned to stay positive. - \textit{Ann. Inst. H. Poincare,
Probab. Statist.}, \textbf{44} (2008), 170--190.

\bibitem{Don12} Doney R.A. Local behavior of first passage probabilities. -
\textit{Probab. Theory Relat. Fields}, \textbf{152} (2012), 559--588.

\bibitem{DGV2004} Dyakonova E.E., Geiger J., Vatutin V.A. On the
survival probability and a functional limit theorem for branching processes
in random environment. - \textit{Markov Process. Relat. Fields}, \textbf{10}
(2004), 289--306.

\bibitem{FE} Feller W.
\newblock{\em An Introduction to Probability Theory
and its Applications.} V.2, Willey, New York-London-Sydney-Toronto, 1971.

\bibitem{GK2000} Geiger J., Kersting G. The survival probability of a
critical branching process in random environment. - \textit{Theory Probab.
Appl.}, \textbf{45} (2000), 607--615.

\bibitem{Koz76} Kozlov M.V. On the asymptotic behavior of the probability of
non-extinction for critical branching processes in a random environment. -
\textit{Theory Probab. Appl.}, \textbf{21} (1976), 791--804.

\bibitem{Koz1995} Kozlov M.V. \ A conditional function limit theorem for a
critical branching process in a random medium. - \textit{Dokl. Akad. Nauk,}
\textbf{344} (1995), 12--15. (In Russian.)


\bibitem{Sen76} Seneta E. \textit{Regularly varying functions.} Lecture
Notes in Mathematics. V.\textbf{508}. Springer, 1976.


\bibitem{Vat2002} Vatutin V.A. Reduced branching processes in random
environment: The critical case. - \textit{Theory Probab. Appl.} \textbf{47}
(2002), 99--113.

\bibitem{Vat2016} Vatutin V. Subcritical branching processes in random
environments. - \textit{Lecture Notes in Statistics - Proceedings Springer, }%
(2016). (In print.)

\bibitem{VD2003} Vatutin V.A., Dyakonova E.E. Galton--Watson branching
processes in random environment, I: Limit theorems. - \textit{Theory Probab.
Appl.}, \textbf{48} (2004), 314--336.

\bibitem{VDS2013} Vatutin V.A., Dyakonova E.E., Sagitov S. Evolution of
branching processes in a random environment. - \textit{Proc. Steklov Inst.
Math.}, \textbf{282} (2013), 220--242.

\bibitem{VL15} Vatutin V., Liu Q. Limit theorems for decomposable
branching processes in random environment. - \textit{J. Appl. Probab.},
\textbf{52} (2015), 877--893.

\bibitem{VW09} Vatutin V.A., Wachtel V. Local probabilities for random
walks conditioned to stay positive. - \textit{Probab. Theory Related Fields,} \textbf{143} (2009), 177--217.

\bibitem{Zol57} Zolotarev V.M. Mellin-Stiltjes transform in probability
theory. - \textit{\ Theory Probab. Appl.,} \textbf{2} (1957), 433--460.
\end{thebibliography}
\end{document}